\documentclass[12pt,twoside,a4paper]{amsart}
\usepackage[centertags]{amsmath}
\usepackage{amssymb,amscd,amsfonts,latexsym,a4wide,amsthm}

\advance\headheight by 3pt

\renewcommand{\subjclass}[1]{\thanks{\emph{2000 Mathematics Subject Classification:}~#1}}
\renewcommand{\keywords}[1]{\thanks{\emph{Keywords and Phrases:}~#1}}

\newtheorem{thm}{Theorem}[section]

\newtheorem{lemma}[thm]{Lemma}

\theoremstyle{definition}

\theoremstyle{remark}

\numberwithin{equation}{section}

\newcommand{\thmref}[1]{Theorem~\ref{#1}}
\newcommand{\secref}[1]{Section~\ref{#1}}
\newcommand{\lemref}[1]{Lemma~\ref{#1}}

\newcommand{\formref}[1]{(\ref{#1})}

\newcommand{\Cc}{{\mathbb C}}

\newcommand{\Zz}{{\mathbb Z}}
\newcommand{\Qq}{{\mathbb Q}}

\newcommand{\Rr}{{\mathbb R}}

\renewcommand{\leq}{\leqslant}
\renewcommand{\geq}{\geqslant}

\newcommand{\kdots}{,\ldots ,}

\renewcommand{\a}{\alpha}
\renewcommand{\b}{\beta}
\newcommand{\ai}{\alpha^{(i)}}
\newcommand{\aj}{\alpha^{(j)}}
\newcommand{\ak}{\alpha^{(k)}}
\newcommand{\ar}{\alpha^{(r)}}
\newcommand{\all}{\{ 1,\ldots ,r\}}

\renewcommand{\proof}{\vskip0.2cm\noindent {\bf Proof.}\,\ }

\title[Distances between the conjugates of an algebraic number]
{Distances between the conjugates of an algebraic number}


\author[J.-H.~EVERTSE]{Jan-Hendrik~EVERTSE}

\begin{document}

\maketitle

\begin{center}
\emph{In memory of Professor B\'{e}la Brindza}
\\[0.5cm]
\end{center}

\begin{abstract}
Let $K$ be a given number field of degree $r\geq 3$,
denote by $\xi\mapsto \xi^{(i)}$ $(i=1\kdots r)$
the isomorphic embeddings
of $K$ into $\Cc$, and let $\Sigma$ be a 
subset of $\all$ of cardinality at least $2$.
Denote by $M(\a )$ the Mahler measure of an algebraic number
$\alpha$. By an elementary argument one shows that
(*) $\prod_{\{ i,j\}\subset\Sigma} |\ai -\aj |
\geq C\cdot M(\a )^{-\kappa }$ 
holds for all $\a$ with $K=\Qq (\a )$, 
with $C=2^{-r(r-1)/2}$ and $\kappa =r-1$.
In the present paper we deduce inequalities (*) with $\kappa <r-1$ and
with a constant $C$ depending on $K$ which are valid for all $\alpha$
with $\Qq (\a )=K$. We obtain such inequalities
with an ineffective constant $C$, using arguments and results
from \cite{Ev1}, \cite{Ev2}, and with an effective constant $C$
using a result from \cite{EvGy}.
 
Define $\kappa (\Sigma )$ to be the infimum of all real numbers $\kappa$
for which there exists a constant $C>0$ such that (*) holds  
for every $\a$ with $\Qq (\a )=K$.
Then clearly $\kappa (\Sigma )\leq r-1$.
We describe the sets $\Sigma$ for which
$\kappa (\Sigma )=r-1$ and we give upper bounds for $\kappa (\Sigma )$
in case that it is smaller than $r-1$. For cubic fields we give the
precise value of $\kappa (\Sigma )$ for each set $\Sigma$. This solves
a problem posed by Mignotte and Payafar \cite[p. 187]{MiPa}.
\end{abstract}

\section{Introduction}\label{1}

Given an algebraic number $\a$ of degree $r$, we denote by
$\a^{(1)}\kdots \ar$ the conjugates of $\a$. 
Letting $a_0$ be the positive integer such that the polynomial
$a_0\prod_{i=1}^r (X-\ai )$ has integer coefficients with greatest common
divisor $1$, we define the Mahler measure and discriminant of $\a$ by
\begin{eqnarray}
\label{1.1}
&M(\alpha ):= 
a_0\prod_{i=1}^r\max\big( 1,|\ai |\big)\, ,& 
\\[0.2cm]
\label{1.2}
&D(\alpha ):=
a_0^{2r-2}
\prod_{1\leq i<j\leq r} 
\big(\ai -\aj\big)^2\, ,&
\end{eqnarray}
respectively.

Let $\Sigma$ be a subset of $\all$ of cardinality $|\Sigma |\geq 2$. 
Then, taking the product over all $2$-element subsets of $\Sigma$,
\begin{eqnarray}\label{1.3}
\prod_{\{ i ,j\}\subset\Sigma} |\ai -\aj |
&\geq&
\prod_{1\leq i<j\leq r} 
\frac{|\ai -\aj |}  
{2\max (1,|\ai |)\max (1,|\aj |)}
\\
\nonumber
&=&2^{-r(r-1)/2}|D(\alpha )|^{1/2}M(\alpha )^{1-r}
\\
\nonumber
&\geq& 
2^{-r(r-1)/2}M(\alpha )^{1-r}
\end{eqnarray}
where the last inequality follows from the fact that $D(\alpha )$
is a non-zero integer.
 
Our purpose is to obtain
improvements of \formref{1.3} with an exponent on $M(\a )$ larger than $1-r$.
More specifically, one could think of improvements
\begin{equation}\label{1.4a}
\prod_{\{ i ,j\}\subset\Sigma} |\ai -\aj |\geq C(r)M(\a)^{-\kappa}
\end{equation}
with $\kappa <r-1$ and a constant $C(r)>0$ depending only on $r$ 
which are
valid for all algebraic numbers of degree $r$,
or, for a given number field $K$ of degree $r$,
\begin{equation}\label{1.4}
\prod_{\{ i ,j\}\subset\Sigma} |\ai -\aj |\geq C(K)M(\a)^{-\kappa}
\end{equation}
with $\kappa <r-1$ and a constant $C(K)>0$ depending on $K$,
which are valid for all $\a$ with $\Qq (\a )=K$.
Apart from a few special cases settled in the literature,
it seems to be difficult to obtain improvements 
of the shape \formref{1.4a}.
In this paper we consider only \formref{1.4}.

We recall some results from the literature dealing with the case
$|\Sigma |=2$,
i.e., inequalities of the shape
\begin{equation}\label{1.3a}
|\ai -\aj |\geq C\cdot M(\alpha )^{-\kappa}\, ,
\end{equation}
where $\Sigma =\{ i,j\}$, $\kappa <r-1$ and either $C=C(r)$ where
$r=\deg\a$ or $C=C(K)$ where $K=\Qq (\a )$.
Mignotte and
Payafar \cite[Theorems 1,2]{MiPa} proved \formref{1.3a} with 
$\kappa =(r-1)/2$ and $C=2^{1-r(r-1)/4}$
if $\ai,\aj\not\in\Rr$ and $\aj\not=\overline{\ai}$;
with $\kappa =(r-1)/3$ and $C=2^{(4-r(r-1))/6}$ 
if $\ai\in\Rr$, $\aj\not\in\Rr$; and with
$\kappa =2$ and $C=2^{1-r}$ if $\Qq (\a )/\Qq$ is a normal extension. Further,
the author \cite[Theorem 4]{Ev1} obtained \formref{1.3a} with 
$\kappa =\frac{41}{42}(r-1)$ and with a constant $C=C(K)$ depending
on $K=\Qq (\alpha )$,
where no restrictions on $\Qq (\a )$, $\ai$, $\aj$ 
are imposed. 
Here $C$ is not effectively computable from the method of proof.
Let $\kappa (r)$ be the infimum of all $\kappa$ for which there is a constant
$C$ such that \formref{1.3a} holds for all algebraic numbers $\alpha$
of degree $r$ and all $i,j$. 
Computations of Collins \cite{Col} suggest that $\kappa (r)=r/2$.
Bugeaud and Mignotte \cite{BuMi} gave an example showing that
if $r$ is even and $r\geq 6$ then $\kappa (r)\geq r/2$.
More generally,
Bugeaud and Mignotte
gave an example showing that for all integers $k,n$ with $k\geq 2$,
$n\geq 3$ 
there are algebraic numbers $\a$ of degree $r=kn$ and of arbitrarily large
Mahler measure, and sets $\Sigma$ of cardinality $k$, such that
\[
\prod_{\{ i ,j\}\subset\Sigma} |\ai -\aj | < c(n,k)M(\a)^{-(1-k^{-1})r}\, .
\]
 
Estimates for the distances between the conjugates of an algebraic number
play an important role in complexity analyses of algorithms for polynomials.
Further, they are of crucial importance in the study of the difference
$w_n(\xi )-w_n^*(\xi )$, 
where $w_n(\xi )$, $w_n^*(\xi )$ are quantities introduced by Mahler
and Koksma, respectively, measuring how well a given transcendental
complex number $\xi$
can be approximated by algebraic numbers of degree $n$,
see the two recent papers by Bugeaud \cite{Bu1},\cite{Bu2}.

In the present paper we are seeking for improvements of 
the shape \formref{1.4}.
Thus, let $K$ be a given number field of degree $r\geq 3$.
Denote by $\xi\mapsto \xi^{(i)}$ ($i=1\kdots r$) the isomorphic embeddings
of $K$ into $\Cc$. The embedding $\xi\mapsto \xi^{(i)}$
is called real if it 
maps $K$ into $\Rr$ and complex if it does not map $K$ into $\Rr$. 
Further, two embeddings $\xi\mapsto \xi^{(i)}$, $\xi\mapsto\xi^{(j)}$ 
are called complex conjugate if 
$\xi^{(j)}=\overline{\xi^{(i)}}$ for $\xi\in K$.  
\\[0.3cm]
{\bf Definition.}
Let $\Sigma$ be a subset of $\all$ of cardinality $\geq 2$.
We define $\kappa (\Sigma )$ to be the infimum of all reals $\kappa$ with the
property that there exists a constant $C(K)>0$ such that
\\[0.2cm]
\formref{1.4}
$\quad
\displaystyle{\prod_{\{ i,j\}\subset\Sigma} |\ai -\aj |
\geq C(K)\cdot M(\a )^{-\kappa }}\quad$
for every $\a$ with $\Qq (\a )=K$.
\\[0.3cm]
From \formref{1.3} it is clear that $\kappa (\Sigma )\leq r-1$.
If $K$ is a cubic field, 
it is possible to give the exact values for
the quantities $\kappa (\Sigma )$.
Our first result is as follows.
\vskip0.5cm

\begin{thm}\label{th:1.1}
Let $K$ be a number field of degree $3$, and $\Sigma$ a subset of $\{ 1,2,3\}$.

(i) Suppose that either $\Sigma =\{ 1,2,3\}$, or $K$ is totally real
and $|\Sigma |=2$, or $\Sigma =\{ i,j\}$ where $\xi\mapsto \xi^{(i)}$
and $\xi\mapsto\xi^{(j)}$ are complex
conjugate. Then $\kappa (\Sigma )=2$.

(ii) Suppose that $\Sigma =\{ i,j\}$, where one of the embeddings
$\xi\mapsto\xi^{(i)}$, $\xi\mapsto\xi^{(j)}$ is real and the other complex. 
Then $\kappa (\Sigma )=\frac{2}{3}$.
\end{thm} 
\vskip0.5cm
We mention that this result solves a problem of Mignotte and Payafar
\cite[bottom of p. 187]{MiPa}.

In the case that the number field $K$ has degree $r\geq 4$, 
we have been able to determine which sets
$\Sigma$ have $\kappa (\Sigma )=r-1$ and to give non-trivial 
(but far from best possible) upper bounds
for $\kappa (\Sigma )$ for the other sets $\Sigma$.
\vskip 0.5cm
 
\begin{thm}\label{th:1.2}
Let $K$ be a number field of degree $r\geq 4$, and $\Sigma$ a subset of 
$\all$.
 
(i) Suppose that either
$\Sigma =\all$ or $\Sigma =\all\backslash\{i_0\}$ where
$\xi\mapsto \xi^{(i_0)}$ is real. Then 
$\kappa (\Sigma )=r-1$.

(ii) Suppose that either $2\leq |\Sigma |\leq r-2$ or
$\Sigma =\all\backslash\{i_0\}$ where
$\xi\mapsto\xi^{(i_0)}$ is complex. Then
\[
\kappa (\Sigma )\leq r-1-\frac{(r-|\Sigma |)^2}{135r}.
\]
\end{thm}
\vskip0.5cm

For instance if $|\Sigma |=2$ part (ii) gives
$\kappa (\Sigma )\leq r-1-(r-2)^2/135r =r-1-O(r)$ which is comparable
to the author's result $\kappa (\Sigma )\leq \frac{41}{42}(r-1)$ mentioned
above. In the other extremal situation $|\Sigma |=r-1$ part (ii) gives
$\kappa (\Sigma )\leq r-1-1/135r$.

Our proof of part (ii) of \thmref{th:1.2} is ineffective.
More precisely,
we prove an inequality of the shape \formref{1.4}
where $\kappa =r-1-(r-|\Sigma |)^2/135r$ and $C(K)$
is not effectively computable
by our method of proof. Below we give an effective version,
but obviously with a value of $\kappa$ much closer to $r-1$.
We denote by $D_K$ the discriminant of a number field $K$.   
\vskip 0.5cm

\begin{thm}\label{th:1.3}
Let $K$ be a number field of degree $r\geq 4$ and let $\Sigma$ be a subset
of $\all$ such that either 
$2\leq |\Sigma |\leq r-2$ or $\Sigma =\all\backslash\{ i_0\}$ where
$\xi\mapsto\xi^{(i_0)}$ is complex. 
Then for every $\alpha$ with $\Qq (\alpha )=K$
we have
\[
\prod_{\{ i ,j\}\subset\Sigma} |\ai -\aj |
\geq C(K)\cdot M(\alpha )^{-\kappa }
\]
with
\begin{equation}\label{1.5}
\kappa = r-1- (c_1r)^{-c_2r^4}|D_K|^{-6r^3},\quad
C(K)=\exp\Big( -(c_3r)^{c_4r^4}|D_K|^{2r^3}\Big)
\end{equation}
where $c_1,c_2,c_3,c_4$ are effectively computable absolute constants.
\end{thm}
\vskip 0.5cm

Our proofs consist of modifications of arguments from \cite{Ev2}. 
We prove \thmref{th:1.1} and part (i) of \thmref{th:1.2} in \secref{2}.
Further, we prove part (ii) of \thmref{th:1.2} and \thmref{th:1.3}
in \secref{3}.

In our proofs we use properties
of equivalence classes of algebraic numbers.
Two algebraic numbers $\a$, $\a^*$ are called equivalent if 
\[
\a^* =\frac{a\a +b}{c\a +d}\quad\text{for some }
\binom{a\ \ b}{c\ \ d}\in{\rm GL}(2,\Zz ).
\]
In \secref{2} we show that if $\Sigma$ satisfies the conditions
of part (i) of \thmref{th:1.2},
then for every $\delta >0$ and every $\a^*$ with $\Qq (\a^* )=K$ there are
infinitely many $\a$ which are equivalent to $\a^*$ and satisfy
\[
\prod_{\{ i,j\}\subset\Sigma} |\ai -\aj |\leq M(\a )^{1-r+\delta}\, .
\]
This implies at once that $\kappa (\Sigma )=r-1$. We use an argument from
\cite{Ev2}, based on Roth's Theorem. The proof of \thmref{th:1.1}
is along the same lines. 

Two equivalent algebraic numbers have the same discriminant.
The author \cite{Ev1} proved that every algebraic number $\a$ with
$\Qq (\a )=K$ is equivalent to an algebraic number $\a^*$
such that
\begin{equation}\label{1.6} 
M(\a^* )\leq A(K)|D(\a )|^{21/(r-1)},
\end{equation}
where $A(K)$ is some ineffective constant depending on $K$.
Thus in \formref{1.3} we may replace the term $|D(\a )|^{1/2}$
by a positive power of $M(\a^* )$,
but $M(\a^* )$ may be much smaller than $M(\a )$.

Provided $\Sigma$ satisfies
the conditions from part (ii)
of \thmref{th:1.2}, we deduce
a refinement of \formref{1.3} 
(\lemref{le:3.3} in \secref{3}) which allows us
to replace
the positive power of $M(\a^*)$ coming from the discriminant by
a positive power of $M(\a )$. This yields at once our upper bound
for $\kappa (\Sigma )$. 
 
To prove \thmref{th:1.3}, we use a result by
Gy\H{o}ry and the author \cite{EvGy}, stating that every algebraic number
$\a$ is equivalent to a number $\a^*$ with 
\begin{equation}\label{1.6a}
M(\a^* )\leq A(K)|D(\a )|^{a(K)}
\end{equation}
where both $A(K)$, $a(K)$ are effectively computable in terms of $K$.
Then the proof of \thmref{th:1.3} is completed similarly as that
of part (ii) of \thmref{th:1.2}.

We mention that both \formref{1.6} and \formref{1.6a} 
were deduced from an inequality of the following type.
Let $K$ be a number field of
degree $r$ and $a,b,c$ non-zero integers of $K$ with $a+b=c$. Then
\begin{equation}\label{1.7}
\prod_{i=1}^r\max (|a^{(i)}|,|b^{(i)}|,|c^{(i)}|)
\leq U\cdot |N_{K/\Qq}(abc)|^{V}\, ,
\end{equation}
where $\xi\mapsto\xi^{(i)}$ $(i=1\kdots r)$ denote as usual the isomorphic
embeddings of $K$ into $\Cc$, and $U,V$ are constants.
Inequality \formref{1.6} follows from a version of \formref{1.7} in which $V=1+\varepsilon$
for any $\varepsilon >0$
and $U=U(K,\varepsilon )$ is some ineffective constant
(see \cite[Lemma 11]{Ev1}). This version
is in turn a consequence of Roth's Theorem over number fields.
Inequality \formref{1.6a} was deduced from a version of \formref{1.7}
in which both $U,V$ are effectively computable in terms
of $K$, but $V$ is rather large 
(see \cite[Theorem]{Gy}, \cite[Corollary]{BuGy}).
The latter is proved by means of linear forms in logarithms
estimates.

As mentioned before, it is as yet open to obtain an inequality of the
shape \formref{1.4a} with $\kappa <r-1$ 
and some constant $C(r)$ depending on $r$. 
We discuss how this is related to certain other open problems.
Assume $\Sigma$ satisfies the condition of part (ii) of \thmref{th:1.2}.
Then by the same reasoning as in the proof of part (ii) of \thmref{th:1.2}
it would be possible to deduce \formref{1.4a} 
with $\kappa =\kappa (r)<r-1$ and $C(r)>0$ 
from an inequality of the shape
\begin{equation}\label{1.6c}
M(\a^* )\leq A(r)|D(\a )|^{a(r)}
\end{equation}
for some
$\a^*$ equivalent to $\a$, where $A(r)$, $a(r)$ depend only on $r$.
Speculating further, by going through the arguments from \cite{Ev1}
it would be possible to deduce \formref{1.6c} 
from a version of \formref{1.7} in which
\[
U=c_1(r)|D_K|^{c_2(r)},\quad V=c_3(r)
\]
where $c_1(r)$, $c_2(r)$, $c_3(r)$ depend only on $r$.
We mention that such a version, with ineffective $c_1(r)$
and effective $c_2(r)$, $c_3(r)$,
can be deduced for instance from a sharpening of Roth's
Theorem over number fields conjectured by Vojta
\cite[\S3, p.65]{Vo}. 
\vskip 1cm

\section{Proofs of \thmref{th:1.1} and part (i) of \thmref{th:1.2}}\label{2}

Our basic tool is the following.
\vskip0.5cm

\begin{lemma}\label{le:2.1}
Let $\a$ be a real, irrational algebraic number and let $\beta_1\kdots\beta_n$
be different complex numbers different from $\a$.
Then for every $\delta >0$ and every $Q$ which is
sufficiently large in terms of $\delta$,
there is a matrix $\binom{a\ b}{c\ d}\in{\rm GL}(2,\Zz )$ such that
\begin{equation}\label{2.1}
\left\{
\begin{array}{rcl}
Q^{-1-\delta}\leq &|\alpha a +b|,\,\, |\alpha c +d|&\leq Q^{-1+\delta}\, ,
\\
Q^{1-\delta}\leq &|\beta_ia+b|,\,\, |\beta_ic+d| &\leq Q^{1+\delta}\quad (i=1\kdots n).
\end{array}\right.
\end{equation}
\end{lemma}

\proof
This lemma is a special case of \cite[Lemma 4.4]{Ev2}. For convenience of the
reader we give the proof.

First we prove the following assertion. For every $\varepsilon$ 
with $0<\varepsilon <1/2$ and every sufficiently
large $Q$, the following holds: if $(x,y)$ is any non-zero point of $\Zz^2$
satisfying
\begin{equation}\label{2.2}
|\a x+y|\leq Q^{-1+\varepsilon},\quad
|\b_i x+y|\leq Q^{1+\varepsilon}\,\,\,(i=1\kdots n),
\end{equation}
then $(x,y)$ satisfies also
\begin{equation}
\label{2.3}
|\a x+y|\geq Q^{-1-2\varepsilon},\quad
|\b_i x+y|\geq Q^{1-2\varepsilon}\,\,\,(i=1\kdots n).
\end{equation}

Below, constants implied by the Vinogradov symbols $\ll$, $\gg$ depend on 
$\alpha$, $\beta_1\kdots\beta_n$ and $\varepsilon$. Let $(x,y)$ be a non-zero
point in $\Zz^2$ satisfying \formref{2.2} but not \formref{2.3}.
Then $x\not= 0$.
First assume that $|\a x+y|< Q^{-1-2\varepsilon}$. Then from \formref{2.2}
we infer $|x|\ll Q^{1+\varepsilon}$ and so
\[
|\a x+y|\ll |x|^{-(1+2\varepsilon )/(1+\varepsilon )}.
\]
By Roth's Theorem, $|x|$ is bounded. But then,
$Q$ is bounded for otherwise 
there are fixed integers $x,y$ with $x\not=0$ satisfying
\formref{2.2} for arbitrarily large $Q$, hence $\a x +y=0$,
which contradicts our assumption that $\a\not\in\Qq$.

Now suppose that $|\b_i x+y|< Q^{1-2\varepsilon}$
for some $i$. Then by using the first inequality in \formref{2.2}
twice, we obtain first $|x|\ll Q^{1-2\varepsilon}$ and then
\[
|\a x+y|\ll |x|^{-(1-\varepsilon )/(1-2\varepsilon )}.
\]
Again by Roth's Theorem, $|x|$ and hence
$Q$ is bounded. This proves 
our assertion.

Now consider the symmetric convex body $S(Q)\subset\Rr^2$, given by
\[  
|\a x+y|\leq Q^{-1},\quad
|\b_i x+y|\leq Q\,\,\,(i=1\kdots n).
\]
$S(Q)$ contains the set of points $(x,y)\in\Rr^2$ with $|\a x +y |\leq Q^{-1}$,
$|y|\ll Q$, therefore its area is $\gg 1$. So by Minkowski's
Theorem, for the successive minima $\lambda_1$, $\lambda_2$ of $S(Q)$
we have
\begin{equation}\label{2.4}
\lambda_1\lambda_2 \ll 1.
\end{equation}
Recall that $\Zz^2$ has a basis $(a,b)$, $(c,d)$ 
(i.e., $\binom{a\ b}{c\ d}\in{\rm GL}(2,\Zz )$) 
such that $(a,b)\in\lambda_1S(Q)$,
$(c,d)\in\lambda_2S(Q)$. Here $\lambda_1$, $\lambda_2$, $(a,b)$, $(c,d)$
depend on $Q$.

Let $0<\varepsilon <1/6$.
Assuming $Q$ is sufficiently large
we have $\lambda_1\geq Q^{-2\varepsilon}$, since otherwise
the point $(a,b)$
would satisfy \formref{2.2} but not
\formref{2.3}, contradicting the assertion proved above.
But then by \formref{2.4} we have
$\lambda_2\ll Q^{2\varepsilon}$,
and hence $\lambda_2\leq Q^{3\varepsilon}$,
assuming that $Q$ is large enough to absorb the constant implied by $\ll$.
This means that
both $(a,b)$, $(c,d)$ satisfy \formref{2.2} with $3\varepsilon$
instead of $\varepsilon$, and then by our assertion they
satisfy also 
\formref{2.3} with $3\varepsilon$ instead of $\varepsilon$,
provided $Q$ is sufficiently large. 
Now choose $\varepsilon <\min (1,\delta)/6$. 
Then $\binom{a\ b}{c\ d}$ satisfies \formref{2.1} 
and our lemma follows.\qed
\vskip0.5cm
\noindent
{\bf Proof of part (i) of Theorems \ref{th:1.1} and \ref{th:1.2}.}
Notice that part (i) of \thmref{th:1.1} is precisely part (i) of
\thmref{th:1.2} with $r=3$. We prove parts (i) of Theorems \ref{th:1.1}
and \ref{th:1.2} simultaneously.
 
Let $K$ be a number field of degree $r\geq 3$.
Without loss of generality we assume that either $\Sigma =\all$
or $\Sigma =\all\backslash\{ 1\}$, where $\xi\mapsto\xi^{(1)}$ is real.
As mentioned in \secref{1}, we pick $\a^*$ with $\Qq (\a^* )=K$ and consider
numbers which are equivalent to $\a^*$. Constants implied by $\ll$, $\gg$
depend on $\a^*$, $K$ and another parameter $\delta$ introduced later.
Let $a_0$ be the integer such that $a_0\prod_{i=1}^r (X-\a^{*(i)})$ has
integer coefficients with greatest common divisor $1$. We use that
for the Mahler measures of the numbers equivalent to $\a^*$ we have
\begin{equation}\label{2.5}
\begin{split}
M\left(\frac{a\a^* +b}{c\a^* +d}\right)
&=a_0\prod_{i=1}^r \max\big( |a\a^{*(i)}+b|,|c\a^{*(i)}+d|\big)
\\
&\quad\text{for }
\binom{a\ \ b}{c\ \ d}\in{\rm GL}(2,\Zz )\, .
\end{split}
\end{equation} 

First suppose that $\Sigma =\all$. We consider numbers $\a_d =(\a^* +d)^{-1}$
with $d\in\Zz$. By \formref{2.5} we have $|d|^r\ll M(\a_d )\ll |d|^r$,
so $M(\a_d )$ tends to $\infty$ with $|d|$.
Moreover, for every $d\in\Zz$ we have
\begin{align*}
\prod_{1\leq i<j\leq r} |\a_d^{(i)} -\a_d^{(j)} | &=
\prod_{1\leq i<j\leq r}\frac{|\a^{*(i)}-\a^{*(j)}|}
{|\a^{*(i)}+d|\cdot |\a^{*(j)}+d|}
\ll |d|^{-r(r-1)}
\\
&\ll M(\a_d )^{1-r}\, .
\end{align*}
Hence $\kappa (\Sigma )=r-1$.

Now assume that 
$\Sigma=\all\backslash\{ 1\}$ where $\xi\mapsto\xi^{(1)}$ is real.
We prove that for every $\delta >0$ there are infinitely many numbers
$\a$ which are equivalent to $\a^*$ and satisfy
\begin{equation}\label{2.6}
\prod_{\{ i,j\}\subset\Sigma} |\ai -\aj |\leq M(\a )^{1-r+\delta}\, .
\end{equation}
This proves $\kappa (\Sigma )=r-1$.

Let $\varepsilon>0$ be a number depending on $\delta$, but much smaller
than $\delta$, which will be specified later. Let $Q>1$. According
to \lemref{le:2.1}, assuming that $Q$ is sufficiently large in terms of $\varepsilon$,
there is a matrix $\binom{a\ b}{c\ d}\in{\rm GL}(2,\Zz )$ such that
\begin{equation}\label{2.7}
\left\{
\begin{array}{rcl}
Q^{-1-\varepsilon}\leq &|\a^{*(1)}a+b|,\,\, |\a^{*(1)}c+d|& \leq Q^{-1+\varepsilon}\, ,
\\
Q^{1-\varepsilon}\leq &|\a^{*(i)}a+b|,\,\, |\a^{*(i)}c+d|& \leq Q^{1+\varepsilon}\quad
(i=2\kdots r).
\end{array}\right.
\end{equation}
Let $\alpha_Q=\frac{a\a^* +b}{c\a^* +d}$; then
$\a_Q$ is equivalent to $\a^*$.  
By \formref{2.5}, \formref{2.7} we have 
\begin{equation}\label{2.8}
Q^{r-2-r\varepsilon}\ll M(\a_Q )\ll Q^{r-2+r\varepsilon}\, ,
\end{equation}
where $a_0$ has been inserted into the constants implied by $\ll$.
Further, by \formref{2.7}, \formref{2.8},
\begin{align*}
\prod_{\{ i,j\}\subset\Sigma} |\a_Q^{(i)}-\a_Q^{(j)}| 
&=\prod_{2\leq i<j\leq r} |\a_Q^{(i)}-\a_Q^{(j)}|
\\ 
&=\prod_{2\leq i<j\leq r} 
\frac{|\a^{*(i)}-\a^{*(j)}|}{|\a^{*(i)}c +d|\cdot |\a^{*(j)}c +d|}
\ll Q^{-(r-1)(r-2)(1-\varepsilon )}
\\
&\ll M(\a_Q)^{-(r-1)(r-2)(1-\varepsilon )/(r-2-r\varepsilon )}.
\end{align*}
Now taking $\varepsilon$ sufficiently small in terms of $\delta$
and then letting $Q\to\infty$ we infer that $\a_Q$ satisfies
\formref{2.6}
and, in view of \formref{2.8}, that $M(\a_Q )\to\infty$. Hence
\formref{2.6} has infinitely many solutions equivalent to $\a^*$.
This completes our proof 
of part (i) of Theorems \ref{th:1.1} and \ref{th:1.2}.
\qed
\vskip0.5cm
\noindent
{\bf Proof of part (ii) of \thmref{th:1.1}.} 
Let $K$ be a cubic field. Without loss of generality we assume that
$\Sigma =\{ 1,2\}$, where $\xi\mapsto \xi^{(1)}$ is real,
$\xi\mapsto \xi^{(2)}$ is complex and $\xi^{(3)}=\overline{\xi^{(2)}}$
for $\xi\in K$.

We recall an argument of Mignotte and Payafar \cite{MiPa}.
Let $\a$ with $\Qq (\a )=K$. Then 
\begin{align*}
&|\a^{(1)}-\a^{(3)}|=|\a^{(1)}-\a^{(2)}|,
\\
&|\a^{(2)}-\a^{(3)}|\leq |\a^{(1)}-\a^{(2)}|+|\a^{(1)}-\a^{(3)}|
=2\cdot |\a^{(1)}-\a^{(2)}|,
\end{align*}
hence
\begin{align*}
|\a^{(1)}-\a^{(2)}|
&\geq
\Big(\frac{1}{2}\prod_{1\leq i<j\leq 3} |\ai -\aj |\Big)^{1/3}
= \Big( \frac{1}{2}a_0^{-2}|D(\a )|^{1/2}\Big)^{1/3}
\\  
&\geq 2^{-1/3}M(\a )^{-2/3}
\end{align*}
where $a_0$ has the meaning from \formref{1.1}, \formref{1.2}.
This proves $\kappa (\Sigma )\leq 2/3$.

To prove the reverse inequality we proceed as in the case $\Sigma =\all$
above. Choose $\a^*$ with $\Qq (\a^* )=K$ and for $d\in\Zz$ define
$\a_d =(\a^* +d)^{-1}$. Then by \formref{2.5} we have 
$|d|^3\ll M(\a_d) \ll |d|^3$ for $d\in\Zz$. Therefore,
$M(\a_d)$ tends to $\infty$ as $|d|\to\infty$.
Moreover,
\[
|\a_d^{(1)}-\a_d^{(2)}|=
\frac{|\a^{*(1)}-\a^{*(2)}|}{|\a^{*(1)}+d|\cdot |\a^{*(2)}+d|}\ll |d|^{-2}
\ll M(\a_d)^{-2/3}.
\]
Hence $\kappa (\Sigma )\geq 2/3$. This completes the proof
of \thmref{th:1.1}.\qed
\vskip 1cm

\section{Proofs of part (ii) of \thmref{th:1.2} and \thmref{th:1.3}}\label{3}

We first state two results of crucial importance for us which are easy
consequences of the literature.
Recall that two equivalent algebraic numbers have the same discriminant. 
\vskip0.5cm

\begin{lemma}\label{le:3.1}
Let $K$ be a number field of degree $r\geq 4$. Then
every $\a$ with $\Qq (\a) =K$ is equivalent to a number $\a^*$ for which
\begin{equation}\label{3.1}
M(\a^* )\leq A_1(K)\cdot |D(\a )|^{21/(r-1)},
\end{equation}
where $A_1(K)$ is a constant depending only on $K$ (which is not effectively
computable from our method of proof).
\end{lemma}
\vskip0.5cm

\begin{lemma}\label{le:3.2}
Let $K$ be a number field of degree $r\geq 4$. Then
every $\a$ with $\Qq (\a )=K$ is equivalent to a number $\a^*$ for which
\begin{equation}\label{3.2}
M(\a^* )\leq A_2(K)\cdot |D(\a )|^{a(K)}
\end{equation}
with 
\begin{equation}\label{3.3}
A_2(K)=\exp\big( (c_5r)^{c_6r^4}|D_K|^{8r^3}\big),\quad
a(K)= (c_7r)^{c_8r^4}|D_K|^{6r^3},
\end{equation}
where $c_5,c_6,c_7,c_8$ are effectively computable absolute constants.
\end{lemma}

\proof
These two lemmata follow from results in the literature
stating that every binary form with integer coefficients and non-zero
discriminant is equivalent to a binary form whose height is bounded
above in terms of the discriminant. 
Given $\alpha$ with $\Qq (\a )=K$, let $a_0$ be the positive integer
such that the binary form $F_{\alpha}(X,Y):=a_0\prod_{i=1}^r (X-\ai Y)$
has integer coefficients with greatest common divisor $1$.
Now \lemref{le:3.1}
follows by applying
\cite[Theorem 1]{Ev1} to $F_{\a}$ and \lemref{le:3.2} by applying  
\cite[Theorem 3']{EvGy} to $F_{\a}$.\qed
\vskip 0.5cm

Our last tool is an improvement of \formref{1.3}. 
\vskip 0.5cm

\begin{lemma}\label{le:3.3}
Let $\a$ be an algebraic number of degree $r\geq 4$.
Let $\a^*$ be equivalent to $\a$ and suppose that $M(\a^* )\leq M(\a )$.
Further, let $\Sigma$ be a subset of $\all$ 
such that either $2\leq |\Sigma |\leq r-2$
or $\Sigma =\all\backslash\{ i_0\}$ where $\a^{(i_0)}\not\in\Rr$.
Then
\begin{eqnarray}
\label{3.4}  
&&\prod_{\{ i,j\}\subset\Sigma}|\ai -\aj |
\\
\nonumber
&&\qquad
\geq 2^{-2r^2}\cdot \frac{|D(\a )|^{1/2}}{M(\a )^{r-1}}
\cdot
\max\left(
1\, ,\,
\frac{|D(\a )|^{1/2}}{M(\a^* )^{r-1}}\cdot
\Big( \frac{M(\a)}{M(\a^*)}\Big)^{4(r-|\Sigma |)^2/9r}\right)\, .
\end{eqnarray}
\end{lemma}

\proof
Write
\[
\a^*=\frac{a\a +b}{c\a +d}\quad\text{with }
\binom{a\ \ b}{c\ \ d}\in{\rm GL}(2,\Zz ).
\]
Define
\[
\varphi_i :=\max (|a\ai +b|, |c\ai +d|),\quad
f_i:= \frac{\max (1, |\ai |)}{\varphi_i}\quad(i=1\kdots r),
\]
and
\[
g_{ij}:= \frac{|\ai -\aj |}{\max (1,|\ai |)\max (1, |\aj |)}
\quad(i,j=1\kdots r).
\]

We first deduce some relations and inequalities for these quantities.
Let $a_0$ be the positive integer such that $a_0\prod_{i=1}^r (X-\ai )$
has integer coefficients with greatest common divisor $1$. Then 
\[
M(\a )=a_0\prod_{i=1}^r\max (1,|\ai |),\quad
M(\a^* )=a_0\prod_{i=1}^r\varphi_i\, ,
\]
hence
\begin{equation}\label{3.5}
f_1\cdots f_r =\frac{M(\a )}{M(\a^* )}.
\end{equation}
It is obvious that
\begin{equation}\label{3.6}
g_{ij}\leq 2\quad\text{for } i,j=1\kdots r\, .
\end{equation}
Further, since $ad-bc=\pm 1$ we have
\[
|\ai -\aj |=|(a\ai +b)(c\aj +d)-(a\aj +b)(c\ai +d)|\leq 2\varphi_i\varphi_j\, ,
\]
hence
\begin{equation}\label{3.7}
g_{ij}f_if_j\leq 2\quad\text{for } i,j=1\kdots r\, .
\end{equation}
From \formref{1.1}, \formref{1.2} it is obvious that
\begin{equation}\label{3.8}
\prod_{1\leq i<j\leq r} g_{ij} =\frac{|D(\a )|^{1/2}}{M(\a )^{r-1}}
\end{equation}
and together with \formref{3.5} this implies
\begin{equation}\label{3.9}
\prod_{1\leq i<j\leq r} (g_{ij}f_if_j) =\frac{|D(\a )|^{1/2}}{M(\a^* )^{r-1}}.  
\end{equation}
Lastly, let $i,j\in\all$ be such that $f_i\leq f_j$.
By \formref{3.5} there is $k\in\all$ with 
$f_k\geq (M(\a )/M(\a^* ))^{1/r}$.
From the vector identity
\[
(\ai -\aj)\Big(\begin{array}{c}1\\ \ak\end{array}\Big)=
(\ai -\ak)\Big(\begin{array}{c}1\\ \aj\end{array}\Big)+
(\ak -\aj)\Big(\begin{array}{c}1\\ \ai\end{array}\Big)
\]
we infer
\begin{align*}
&|\ai-\aj |\cdot\max (1,|\ak |)
\\
&\qquad\leq
|\ai-\ak |\cdot\max (1,|\aj |)+|\ak-\aj |\cdot\max (1,|\ai |)
\end{align*}
and so $g_{ij}\leq g_{ik}+g_{kj}$. Now invoking \formref{3.7} 
and our assumption $f_i\leq f_j$ we obtain
\[
g_{ij}f_if_jf_k\leq g_{ik}f_if_kf_j+g_{kj}f_kf_jf_i\leq 2f_j+2f_i\leq 4f_j,
\]
and by dividing by $f_j$ and using our assumption on $k$ we arrive at 
\begin{equation}\label{3.10}
g_{ij}f_i\cdot \Big(\frac{M(\a )}{M(\a^*)}\Big)^{1/r}\leq 4\,\,\,
\mbox{for $i,j\in\all$ with $f_i\leq f_j$.}
\end{equation} 

Having finished our preparations, we now commence with our proof.
By \formref{3.8} we have
\[
\prod_{\{ i,j\}\subset\Sigma} |\ai -\aj | 
\geq \prod_{\{ i,j\}\subset\Sigma} g_{ij}
=\frac{|D(\a )|^{1/2}}{M(\a )^{r-1}}
\cdot \prod_{\{ i,j\}\not\subset\Sigma} g_{ij}^{-1}\, .
\]
By \formref{3.6} we have  
$\prod_{\{ i,j\}\not\subset\Sigma} g_{ij}^{-1}\geq 2^{-r(r-1)/2}$.
So in order to prove \formref{3.4}, it suffices to prove that
\begin{equation}\label{3.11}
\prod_{\{ i,j\}\not\subset\Sigma} g_{ij}^{-1}\geq
2^{-2r^2}\frac{|D(\a )|^{1/2}}{M(\a^* )^{r-1}}\cdot
\left(\frac{M(\a )}{M(\a^* )}\right)^{4(r-|\Sigma |)^2/9r}\, .
\end{equation}
We distinguish two cases.

First assume that $2\leq |\Sigma |\leq r-2$. Put $l:= r-|\Sigma |$.
Choose $j_0\in \Sigma$. Without loss of generality we may
assume that $\{ j_0\}\cup \all\backslash\Sigma =\{ 1\kdots l+1\}$ and that
\begin{equation}\label{3.12}  
f_1\leq f_2\leq\cdots \leq f_{l+1}.
\end{equation}
Notice that if $1\leq i<j\leq l+1$ then $\{ i,j\}\not\subset\Sigma$.
Denote by $A$ the collection of pairs of indices
$(i,j)$ with $2\leq i<j\leq \min (2i-1,l+1)$ and by $B$ the collection
of pairs $(i,j)$ such that $1\leq i<j\leq r$, $(i,j)\not\in A$ and 
$\{ i,j\}\not\subset\Sigma$.
By an easy computation we have $|A|= l^2/4$ if $l$ is even,
$|A|=(l^2-1)/4$ if $l$ is odd and so for both $l$ even or odd
(using $l\geq 3$ if $l$ is odd),
\begin{equation}\label{3.13}
|A|\geq 2l^2/9 =2(r-|\Sigma |)^2/9\, .
\end{equation}
First take $(i,j)\in A$. 
Then $1\leq 2i-j<i<j\leq l+1$,
and so by \formref{3.10}, \formref{3.12},
\begin{align*}
g_{ij}^{-1}&\geq g_{ij}^{-1}\cdot 
\frac{1}{4}g_{2i-j,i}f_{2i-j}\Big(\frac{M(\a )}{M(\a^*)}\Big)^{1/r}
\cdot
\frac{1}{4}g_{ij}f_i\Big(\frac{M(\a )}{M(\a^*)}\Big)^{1/r}
\\
&=\frac{1}{16}g_{2i-j,i}f_{2i-j}f_i\cdot \Big(\frac{M(\a )}{M(\a^*)}\Big)^{2/r}.
\end{align*}
For $(i,j)\in B$ we use \formref{3.6}. Thus we obtain
\begin{equation}
\label{3.14}
\prod_{\{ i,j\}\not\subset\Sigma} g_{ij}^{-1}
\geq 2^{-|B|-4|A|}\cdot \Big(\frac{M(\a )}{M(\a^*)}\Big)^{2|A|/r}
\cdot\prod_{(i,j)\in A} (g_{2i-j,i}f_{2i-j}f_i).
\end{equation}
Thanks to the fact that
the sets $\{2i-j,i\}$ $((i,j)\in A)$ are
distinct (which is crucial and the main motivation for our set-up),
we infer from \formref{3.9}, \formref{3.7},
\[
\prod_{(i,j)\in A} (g_{2i-j,i}f_{2i-j}f_i)
\geq 2^{|A|-r(r-1)/2}\cdot \frac{|D(\a )|^{1/2}}{M(\a^* )^{r-1}}.
\]
By inserting this and \formref{3.13} into \formref{3.14},
and using our assumption $M(\a )\geq M(\a^* )$ we arrive at
\begin{align*}
\prod_{\{ i,j\}\not\subset\Sigma} g_{ij}^{-1}
&\geq 2^{-|B|-3|A|-r(r-1)/2}\cdot \frac{|D(\a )|^{1/2}}{M(\a^* )^{r-1}}
\cdot \Big(\frac{M(\a )}{M(\a^*)}\Big)^{2|A|/r}
\\
&\geq 2^{-2r^2}\cdot\frac{|D(\a )|^{1/2}}{M(\a^* )^{r-1}}\cdot
\Big(\frac{M(\a )}{M(\a^*)}\Big)^{4(r-|\Sigma |)^2/9r}
\end{align*}
which is \formref{3.11}.

We now treat the case $\Sigma =\all\backslash\{ i_0\}$ 
where $\alpha^{(i_0)}\not\in\Rr$. Without loss of generality we assume that
$i_0=1$ and that $\alpha^{(2)}=\overline{\alpha^{(1)}}$. Then $f_1=f_2$
and so by \formref{3.10},
\begin{align*}
g_{12}^{-1}&\geq g_{12}^{-1}\cdot 
\frac{1}{4}g_{12}f_1\Big(\frac{M(\a )}{M(\a^*)}\Big)^{1/r}
\cdot
\frac{1}{4}g_{12}f_2\Big(\frac{M(\a )}{M(\a^*)}\Big)^{1/r}
\\
&=\frac{1}{16}\cdot \Big(\frac{M(\a )}{M(\a^*)}\Big)^{2/r} 
\cdot g_{12}f_1f_2.
\end{align*}
Now by \formref{3.9}, \formref{3.7} we have
\[
g_{12}f_1f_2\geq 2^{1-r(r-1)/2}\frac{|D(\a )|^{1/2}}{M(\a^* )^{r-1}}.
\]
Hence
\begin{align*}
\prod_{\{ i,j\}\not\subset\Sigma} g_{ij}^{-1}
&=\prod_{j=2}^r g_{1j}^{-1}
\geq 2^{2-r}g_{12}^{-1}\geq 2^{-2-r}g_{12}f_1f_2
\cdot \Big(\frac{M(\a )}{M(\a^*)}\Big)^{2/r}
\\
&\geq 2^{-1-r-r(r-1)/2}\cdot\frac{|D(\a )|^{1/2}}{M(\a^* )^{r-1}}\cdot
\Big(\frac{M(\a )}{M(\a^*)}\Big)^{2/r}
\end{align*}
which implies \formref{3.11}. This completes the
proof of \lemref{le:3.3}.\qed 
\vskip0.5cm

In what follows, Let $K,\Sigma$, $r$ be as in part (ii) of \thmref{th:1.2}.
Take $\a$ with $\Qq (\a )=K$. From the equivalence
class of $\a$ we choose an element $\a^*$
of minimal Mahler measure. Thus, $M(\a^*)\leq M(\a )$ hence
all conditions of \lemref{le:3.3} are satisfied.
Further, $\a^*$ satisfies the inequalities \formref{3.1} and
\formref{3.2} in \lemref{le:3.1}, \lemref{le:3.2}, respectively.
Put 
\[
u:= 4(r-|\Sigma |)^2/9r.
\]
Let $0\leq\theta \leq 1$. Then \formref{3.4} implies
\begin{eqnarray}\label{3.15}
&&\prod_{\{ i,j\}\not\in\Sigma} |\ai -\aj |
\\
\nonumber
&&\qquad
\geq
2^{-2r^2}\frac{|D(\a )|^{1/2}}{M(\a )^{r-1}}\cdot
\left(\frac{|D(\a )|^{1/2}}{M(\a^* )^{r-1}}M(\a )^uM(\a^* )^{-u}\right)^{\theta}
\\
\nonumber 
&&\qquad =
2^{-2r^2}\cdot |D(\a )|^{(1+\theta )/2}\cdot M(\a^*)^{-\theta (r-1+u)}
\cdot M(\a )^{1-r+\theta u}\, .
\end{eqnarray}
We prove part (ii) of \thmref{th:1.2} and \thmref{th:1.3} by
combining \formref{3.15} with
\formref{3.1}, \formref{3.2}, respectively,
and choosing an appropriate value for $\theta$.
\vskip0.5cm
\noindent  
{\bf Proof of part (ii) of \thmref{th:1.2}.}
By \formref{3.1} we have 
\[
|D(\a )|\geq A_1(K)^{-(r-1)/21}M(\a^*)^{(r-1)/21}.
\]
We insert this into \formref{3.15} 
and then choose $\theta$ to make the exponent on $M(\a^*)$ equal to $0$.
Thus,
\begin{align*}
&\prod_{\{ i,j\}\not\in\Sigma} |\ai -\aj |
\\
&\qquad\geq 2^{-2r^2}A_1(K)^{-\frac{r-1}{42}(1+\theta )}\cdot 
M(\a^*)^{\frac{r-1}{42}(1+\theta ) -\theta (r-1+u)}\cdot 
M(\a )^{1-r+\theta u}
\\
&\qquad=2^{-2r^2}A_1(K)^{-\frac{r-1}{42}(1+\theta )}\cdot M(\a )^{1-r+\theta u}\, ,
\end{align*}
where $\frac{r-1}{42}(1+\theta ) =\theta (r-1+u)$, that is,
\[
\theta =\frac{1}{41 +42u/(r-1)}\, .
\]
Consequently, using $u\leq 4(r-1)/9$,
\begin{align*}
\kappa (\Sigma )
&\leq r-1-\theta u \leq r-1-\frac{4(r-|\Sigma |)^2/9r}{41+42\times 4/9}
\\
&\leq
r-1-\frac{(r-|\Sigma |)^2}{135r}\, .
\end{align*}
This proves part (ii) of \thmref{th:1.2}.\qed
\vskip0.5cm
\noindent
{\bf Proof of \thmref{th:1.3}.}
By \formref{3.2} we have 
\[
|D(\a )|\geq A_2(K)^{-1/a(K)}M(\a ^*)^{1/a(K)}.
\]
Similarly as above, we insert this into \formref{3.15}, and choose $\theta$
such that the exponent on $M(\a^* )$ becomes $0$. Thus,
\begin{eqnarray}
\label{3.16}
&&\prod_{\{ i,j\}\not\in\Sigma} |\ai -\aj |
\\
\nonumber
&&\qquad\geq 2^{-2r^2}A_2(K)^{-\frac{1+\theta }{2a(K)}}\cdot 
M(\a^*)^{\frac{1+\theta}{2a(K)} -\theta (r-1+u)}\cdot 
M(\a )^{1-r+\theta u}
\\
\nonumber
&&\qquad = 2^{-2r^2}A_2(K)^{-\frac{1+\theta }{2a(K)}}
\cdot M(\a )^{1-r+\theta u}\, ,
\end{eqnarray} 
where $\frac{1+\theta}{2a(K)} =\theta (r-1+u)$, that is,
\[
\theta =\frac{1}{(2r-2+2u)a(K)-1}.
\]
With this choice of $\theta$ we have
\begin{align*}
2^{-2r^2}A_2(K)^{-(1+\theta)/2a(K)}
&\geq 2^{-2r^2}A_2(K)^{-1/a(K)}
\\
&\geq 2^{-2r^2}
\exp\Big( -(c_5r)^{c_6r^4}|D_K|^{8r^3}(c_7r)^{-c_8r^4}|D_K|^{-6r^3}\Big)
\\
&\geq \exp\Big( -(c_3r)^{c_4r^4}|D_K|^{2r^3}\Big)
\end{align*}
and, using $4/9r\leq u\leq 4(r-1)/9$,
\begin{align*}
\theta u &\geq \frac{4}{9r}
\cdot
\big\{ (2r-2 +\mbox{$\frac{8}{9}$}(r-1))(c_5r)^{c_6r^4}|D_K|^{6r^3}\big\}^{-1}
\\
&\geq (c_1r)^{-c_2r^4}|D_K|^{-6r^3}.
\end{align*}
By inserting this into \formref{3.16}, \thmref{th:1.3} follows.\qed

\vskip0.5cm
{\scshape\small
J.-H. Evertse

Universiteit Leiden, Mathematisch Instituut

Postbus 9512, 2300 RA Leiden, The Netherlands
}

{\itshape\small E-mail address:} {\ttfamily\small evertse@math.leidenuniv.nl}

\end{document}